# Stochastic differential equations with jumps*,†

**Richard F. Bass**

*Department of Mathematics, University of Connecticut*
*Storrs, CT 06269-3009*
*e-mail:* bass@math.uconn.edu

**Abstract:** This paper is a survey of uniqueness results for stochastic differential equations with jumps and regularity results for the corresponding harmonic functions.

## 1. Introduction

Researchers have increasingly been studying models from economics and from the natural sciences where the underlying randomness contains jumps. To give an example from financial mathematics, the classical model for a stock price is that of a geometric Brownian motion. However, wars, decisions of the Federal Reserve and other central banks, and other news can cause the stock price to make a sudden shift. To model this, one would like to represent the stock price by a process that has jumps.

This paper is a survey of some aspects of stochastic differential equations (SDEs) with jumps. As will quickly become apparent, the theory of SDEs with jumps is nowhere near as well developed as the theory of continuous SDEs, and in fact, is still in a relatively primitive stage. In my opinion, the field is both fertile and important. To encourage readers to undertake research in this area, I have mentioned some open problems.

Section 2 is a description of stochastic integration when there are jumps. Section 3 describes a number of the most important types of SDEs that involve jumps. In Section 4 I discuss a few pathwise uniqueness results, and in Section 5 I discuss some results relating to martingale problems. Section 6 concerns the regularity of harmonic functions with respect to the corresponding integral operators.

A forthcoming book by Applebaum [Ap] contains much information on the subject of this paper, as well as applications to finance.

*This is an original survey paper.
†Research partially supported by NSF grant DMS-0244737.

## 2. Stochastic calculus

Before we can describe SDEs with respect to jump processes, we need to talk a bit about the differences between the stochastic calculus for continuous processes and for processes with jumps. Some good references for this are the volumes by Dellacherie and Meyer [DM1], [DM2], Meyer's course [Me], and the books by He, Wang, and Yan [HWY], Jacod [Jd], Protter [Pr], and von Weizsäcker and Winkler [WW]. I have some notes on the web at [Ba5]; my notes [Ba6] provide further background on predictability.

Suppose we are given a probability space $(\Omega, \mathcal{F}, \mathbb{P})$ and a filtration $\{\mathcal{F}_t\}$. We suppose that the filtration satisfies the "usual conditions," which means that $\mathcal{F}_t$ is right continuous and each $\mathcal{F}_t$ is complete. We say that $X_t$ is right continuous with left limits if there exists a null set $N$ so that if $\omega \notin N$, then $\lim_{u\downarrow t} X_u(\omega) = X_t(\omega)$ for all $t$ and $\lim_{s\uparrow t} X_s(\omega)$ exists for all $t$. The French abbreviation for this is "cádlág." Given a process $X_t$ that is right continuous with left limits, let $X_{t-} = \lim_{s\uparrow t} X_s$ and $\Delta X_t = X_t - X_{t-}$.

A stochastic process $X$ can be viewed as a map from $[0,\infty) \times \Omega$ to (usually) $\mathbb{R}$. We define a $\sigma$-field $\mathcal{P}$ on $[0,\infty) \times \Omega$, called the predictable $\sigma$-field or previsible $\sigma$-field, by letting $\mathcal{P}$ be the $\sigma$-field on $[0,\infty) \times \Omega$ generated by the class of stochastic processes that are adapted to the filtration $\{\mathcal{F}_t\}$ and have left continuous paths. One can show that $\mathcal{P}$ is also generated by the class of processes $H(s,\omega) = G(\omega) 1_{(a,b]}(s)$, where $G$ is bounded and $\mathcal{F}_a$-measurable.

Suppose $A_t$ is an $\mathcal{F}_t$-adapted process that is right continuous with left limits whose paths are increasing, and let us suppose for simplicity that $A_t$ is integrable for each $t$. Since $A_t \geq A_s$, then $\mathbb{E}[A_t \mid \mathcal{F}_s] \geq \mathbb{E}[A_s \mid \mathcal{F}_s] = A_s$, so trivially $A_t$ is a submartingale. By the Doob-Meyer decomposition, there exists a predictable increasing adapted process $\widetilde{A}_t$ such that

$$A_t = \text{ martingale } + \widetilde{A}_t.$$

We call $\widetilde{A}_t$ the compensator or dual predictable projection of $A_t$. If $A_t$ has paths of bounded variation, we write $A_t = A_t^+ - A_t^-$, where $A_t^+$ and $A_t^-$ are increasing processes, and then define $\widetilde{A}_t = \widetilde{A_t^+} - \widetilde{A_t^-}$.

Before defining the stochastic integral, we discuss the decomposition of square integrable martingales. A martingale $M_t$ is square integrable if $\sup_t \mathbb{E}\, M_t^2 < \infty$.

If $M_t$ is a square integrable martingale that is right continuous with left limits, for each integer $i$ let $T_{i1} = \inf\{t : |\Delta M_t| \in [2^i, 2^{i+1})\}$, $T_{i2} = \inf\{t > T_{i1} : |\Delta M_t| \in [2^i, 2^{i+1})\}$, and so on. Since $M_t$ is right continuous with left limits, $T_{ij} \to \infty$ as $j \to \infty$ for each $i$. Let $\{S_m\}$ be some ordering of the $\{T_{ij}\}$. We then have a countable sequence of stopping times $S_i$ such that every jump of $M$ occurs at one of the $S_i$ and $\Delta M_{S_i}$ is a bounded random variable for each $i$.

Let us set $A_i(t) = \Delta M_{S_i} 1_{(t \geq S_i)}$ and set $M_i(t) = A_i(t) - \widetilde{A}_i(t)$. One can then prove the following theorem; see [Me], T.II.11 or [Ba5], Th. 2.3.

**Theorem 2.1** *Suppose $M_t$ is a square integrable martingale that is right continuous with left limits and $M_i(t)$ is defined as above. Then each $M_i(t)$ is a* *square integrable martingale and* $\sum_{i=1}^{\infty} M_i(t)$ *converges in* $L^2$ *for each* $t$*. If* $M_t^c = M_t - \sum_{i=1}^{\infty} M_i(t)$*, it is possible to find a version of* $M_t^c$ *that is a square integrable martingale with continuous paths. Moreover,* $M^c$ *and the* $M_i$ *are mutually orthogonal.*

The bit about finding a version of $M^c$ is due to the fact that the infinite sum converges in $L^2$, but there is a null set that depends on $t$. In fact, it is possible to arrange matters so that there a single null set. Saying two martingales $N_1(t)$ and $N_2(t)$ are orthogonal means here that $\mathbb{E}\,[N_1(T)N_2(S)] = 0$ for every pair of stopping times $S$ and $T$.

Recall that if $M_t$ is a martingale, then $\langle M \rangle_t$ is the unique increasing predictable process such that $M_t^2 - \langle M \rangle_t$ is a martingale. The existence and uniqueness of $\langle M \rangle_t$ is guaranteed by the Doob-Meyer decomposition. If $M_t$ is a square integrable martingale whose paths are left continuous with right limits, define

$$[M]_t = \langle M^c \rangle_t + \sum_{s \le t} |\Delta M_s|^2.$$

Here $M^c$ is the continuous part of $M$ given in Theorem 2.1. One can show, using Theorem 2.1, that $M_t^2 - [M]_t$ is a martingale, and in particular, $\mathbb{E}\,M_t^2 = \mathbb{E}\,[M]_t$.

With this as background, we can now proceed to a definition of stochastic integrals with respect to a square integrable martingale. We want our integrands to be predictable. Let us take a moment to explain why this is very natural. Let $M_t = 1 + P_{t\wedge 1} - t \wedge 1$, i.e., a Poisson process minus its mean stopped at time 1, and then with 1 added so that $M_t$ is nonnegative. Let us suppose that $M_t$ is the price of a stock, and $H_s$ is the number of shares we hold at time $s$. With this investment strategy, it is not hard to see that the net profit (or loss) at time 1 is $\int_0^1 H_s dM_s$. Let $T$ be the first time the Poisson process jumps one. If we were allowed to choose $H_s$ to be zero for $s < T$ and 1 for $T \le s$, our profit at time 1 would be 1 on the event that $(T \le 1)$ and 0 on the event $(T > 1)$; we would have made a profit without any risk (if $M_t$ has paths of bounded variation, the stochastic integral and Lebesgue-Stieltjes integral will coincide). The problem is that we looked an instantaneous amount into the future to see when the Poisson process jumped. We can't allow that, and the way to prevent this is to require $H_s$ to be predictable.

Let us give the definition of stochastic integral. Suppose $M_t$ is a square integrable martingale with paths that are right continuous with left limits. If

$$H_s(\omega) = \sum_{i=1}^{m} G_i(\omega) 1_{(a_i,b_i]}(s), \tag{2.1}$$

where $G_i$ is bounded and $\mathcal{F}_{a_i}$-measurable, define

$$N_t = \int_0^t H_s dM_s = \sum_{i=1}^{m} G_i (M_{t\wedge b_i} - M_{t\wedge a_i}).$$

Just as in [Ba3], pp. 43–44, the left hand side will be a martingale and very similarly to the proof of [Ba3], Section I.5 with [ ] instead of $\langle\ \rangle$, $N_t^2 - [N]_t$ will be a martingale, where $[N]_t = \int_0^t H_s^2 d[M]_s$.

If $H$ is only $\mathcal{P}$-measurable and satisfies $\mathbb{E} \int_0^\infty H_s^2 d[M]_s < \infty$, approximate $H$ by integrands $H_s^n$ where each $H_s^n$ is of the form given in (2.1). Define $N_t^n = \int_0^t H_s^n dM_s$. By the same proof as in [Ba3], Section I.5, the martingales $N_t^n$ converge in $L^2$. We call the limit $N_t = \int_0^t H_s dM_s$. The stochastic integral is a square integrable martingale, its paths are right continuous with left limits, its definition is independent of which sequence $H_s^n$ we choose, and $[N]_t = \int_0^t H_s^2 d[M]_s$.

We want to generalize the definition of stochastic integral to more general processes. For example, a Cauchy process, even if stopped at a fixed time, is not square integrable, but we should be able to handle stochastic integrals with respect to a Cauchy process by looking at the large jumps separately.

We say $M_t$ is a local martingale if there exist stopping times $T_n$ increasing to infinity such that for each $n$ the process $M_{t\wedge T_n}$ is a uniformly integrable martingale. A semimartingale is a process of the form $X_t = X_0 + M_t + A_t$, where $X_0$ is a finite random variable that is $\mathcal{F}_0$-measurable, $M_t$ is a local martingale, and $A_t$ is a process whose paths have bounded variation on $[0,t]$ for each $t$. A key result concerning semimartingales is the following reduction theorem; see [Me], T.IV.8 or [Ba5], Th. 5.4. This is easy in the case of continuous semimartingales, but not at all in the case of semimartingales that have jumps.

**Theorem 2.2** *Suppose $X_t$ is a semimartingale. There exist stopping times $S_n \uparrow \infty$ such that $X_{t\wedge S_n} = U_t^n + V_t^n$, where $U^n$ is a square integrable martingale and $V^n$ is a process whose paths have bounded variation and the total variation of $V^n$ over the time interval $[0,\infty]$ is finite. Moreover, $U_t^n = U_{S_n}^n$ and $V_t^n = V_{S_n}^n$ for $t \geq S_n$.*

If $X_t$ is a local martingale and $S_n$ are stopping times such as in Theorem 2.2, set $X_{t\wedge S_n}^c = (U^n)_t^c$ for each $n$ and $[X]_{t\wedge S_n} = \langle X^c\rangle_{t\wedge S_n} + \sum_{s\leq t\wedge S_n} \Delta X_s^2$, where $(U^n)^c$ is defined as in Theorem 2.1. It is not hard to see that these definition are independent of the choice of stopping times $S_n$.

Next we weaken the assumptions on $H$. A predictable process $H_s$ is locally bounded if there exist stopping times $R_n \uparrow \infty$ and constants $K_n$ such that the process $H$ is bounded by $K_n$ on $[0, R_n]$. If $H$ is locally bounded and $X$ is a semimartingale, we define $\int_0^t H_s dX_s$ by setting

$$\int_0^t H_s dX_s = \int_0^t H_{s\wedge R_n} dU_{s\wedge S_n}^n + \int_0^t H_{s\wedge R_n} dV_{s\wedge S_n}^n$$

when $t \leq R_n \wedge S_n$; the first integral is a stochastic integral with respect to the square integrable martingale $U^n$ and the second is a Lebesgue-Stieltjes integral. Since $R_n \wedge S_n \uparrow \infty$, this defines $\int_0^t H_s dX_s$ for all $t$. It can be shown that this integral does not depend on the choice of $R_n$ and $S_n$.

Itô's formula (see [Me], T.III.3 or [Ba5], Th. 4.1) is

**Theorem 2.3** *Suppose $X$ is a semimartingale and $f$ is a $C^2$ function. Then $f(X_t)$ is also a semimartingale and we have*

$$\begin{aligned} f(X_t) &= f(X_0) + \int_0^t f'(X_{s-})dX_s + \tfrac{1}{2}\int_0^t f''(X_{s-})d\langle X^c\rangle_s \\ &\quad + \sum_{s\le t}[f(X_s) - f(X_{s-}) - f'(X_{s-})\Delta X_s]. \end{aligned}$$

Note that $f'(X_{s-})$ is a left continuous process, hence predictable. For $d$-dimensional processes, each of whose components is a semimartingale, the obvious generalization holds.

## 3. Jump processes

*SDEs with respect to Lévy processes.* The classical SDE with respect to Brownian motion is

$$dX_t = \sigma(X_t)dW_t + b(X_t)dt, \tag{3.1}$$

where $W_t$ is a Brownian motion. The simplest analogue of this in the jump case is

$$dX_t = a(X_{t-})dZ_t, \tag{3.2}$$

where $Z_t$ is a Lévy process. To be even more specific, we might take $Z_t$ to be a symmetric stable process of index $\alpha$. Already, even in this very special case, there are several interesting things one can say. Note that we write $X_{t-}$ instead of $X_t$ in order that the integrand be predictable.

It would be natural to allow a process to have both a continuous component and a jump component, so one might want to consider the SDE

$$dX_t = \sigma(X_t)dW_t + b(X_t)dt + a(X_{t-})dZ_t. \tag{3.3}$$

*Poisson point processes.* For many applications, (3.3) is not a very good model. For example, suppose one wants to model a stock price in such a way that the underlying randomness is given by a jump process. For simplicity let us consider (3.2) instead of (3.3). If $Z_t$ has a jump of size $z$, then $X_t$ will have a jump of size $a(X_{t-})z$. However, one might very well want $X_t$ to have a jump whose size depends on $X_{t-}$ and $z$, but is not necessarily linear in $z$. If the underlying randomness has a big jump, the behavior of $X$ might be qualitatively different from the case where the underlying randomness has a small jump.

To obtain models with this extra versatility, we need to consider Poisson point processes. See [Sk] for details. Let $(S,\lambda)$ be an arbitrary measure space (letting $S=\mathbb{R}$ and $\lambda$ be Lebesgue measure will usually do). For each $\omega\in\Omega$ let $\mu(\omega,dt,dz)$ be a measure on $[0,\infty)\times S$. The random measure $\mu$ is a Poisson point process if (i) for each set $A\subset S$ with $\lambda(S)<\infty$ the process $\mu([0,t]\times A)$ is a Poisson process with parameter $\lambda(A)$ and (ii) if $A_1,\ldots,A_n$ are disjoint subsets of $S$ with $\lambda(A_i)<\infty$, then the processes $\mu([0,t]\times A_i)$ are independent.

Define a non-random measure $\nu$ by $\nu([0,t]\times A)=t\lambda(A)$. If $\lambda(A)<\infty$, then $\mu([0,t]\times A)-\nu([0,t]\times A)$ is the same as a Poisson process minus its mean, hence is a martingale.

We can define a stochastic integral with respect to the compensated point process $\mu - \nu$ as follows. Suppose $H(s,z) = H(\omega, s, z)$ is of the form

$$H(s,z) = \sum_{i=1}^n K_i(\omega) 1_{(a_i,b_i]}(s) 1_{A_i}(z), \tag{3.4}$$

where for each $i$ the random variable $K_i$ is $\mathcal{F}_{a_i}$-measurable and $A_i \subset S$ with $\lambda(A_i) < \infty$. We define

$$N_t = \int_0^t \int H(s,z)\,(\mu - \nu)(ds,dz) = \sum_{i=1}^n K_i (\mu - \nu)(((a_i, b_i] \cap [0,t]) \times A_i).$$

By linearity it is easy to see that $N_t$ is a martingale. It is also easy to see that $N^c = 0$ and

$$[N]_t = \int_0^t \int H(s,z)^2 \mu(ds,dz). \tag{3.5}$$

With a little work one can show

$$\langle N \rangle_t = \int_0^t \int H(s,z)^2 \nu(ds,dz). \tag{3.6}$$

Suppose $H(s,z)$ is a predictable process in the following sense: $H$ is measurable with respect to the $\sigma$-field generated by all processes of the form (3.4). Suppose also that $\mathbb{E} \int_0^\infty \int H(s,z)^2 \nu(ds,dz) < \infty$. If we take processes $H^n$ of the form (3.4) converging to $H$ in an appropriate way, the corresponding $N_t^n = \int_0^t \int H^n(s,z) d(\mu - \nu)$ will converge in $L^2$, and we call the limit $N_t$ the stochastic integral of $H$ with respect to $\mu - \nu$. One can show that (3.5) and (3.6) are still valid. One may think of the stochastic integral as follows: if $\mu$ assigns mass one to the point $(t,z)$, then $N_t$ jumps at this time $t$ and the size of the jump is $H(t,z)$.

Now consider a stochastic differential equation with respect to a compensated Poisson point process. Look at

$$dX_t = \sigma(X_t)\,dW_t + b(X_t)\,dt + \int F(X_{t-}, z)\,d(\mu - \nu), \qquad X_0 = x_0. \tag{3.7}$$

This means

$$X_t = x_0 + \int_0^t \sigma(X_s) dW_s + \int_0^t b(X_s) ds + \int_0^t \int F(X_{s-}, z)(\mu - \nu)(ds,dz),$$

where $W_t$ is a standard Brownian motion on $\mathbb{R}$. This is a quite general formulation, as is shown in Çinlar and Jacod [CJ].

*Martingale problems.* For simplicity we consider here the SDE (3.7) with $\sigma \equiv b \equiv 0$. Suppose $f \in C^2$ and suppose that $X_t$ is the solution to

$$dX_t = \int F(X_{t-}, z)\,d(\mu - \nu). \tag{3.8}$$

By Itô's formula (Theorem 2.3),

$$f(X_t) = f(X_0) + \int_0^t f'(X_{s-})dX_s + \sum_{s\le t}[f(X_s) - f(X_{s-}) - f'(X_{s-})\Delta X_s].$$

The stochastic integral term is a martingale. We can write $f(X_s)$ as $f(X_{s-} + \Delta X_s)$. The jump of $X_t$ at time $s$ is equal to $F(X_{s-}, z)$ if $\mu$ assigns mass one to the point $(s, z)$, and so $f(X_t) - f(X_0)$ is equal to a martingale plus

$$\int_0^t \int [f(X_{s-} + F(X_{s-}, z)) - f(X_{s-}) - f'(X_{s-})F(X_{s-}, z)]\mu(ds, dz).$$

This in turn is equal to a martingale plus

$$\begin{aligned}\int_0^t \int &[f(X_{s-} + F(X_{s-}, z)) - f(X_{s-}) - f'(X_{s-})F(X_{s-}, z)]\nu(dz, ds)\\ &= \int_0^t \int [f(X_{s-} + F(X_{s-}, z)) - f(X_{s-}) - f'(X_{s-})F(X_{s-}, z)]\lambda(dz)ds.\end{aligned}$$

If we now set

$$\mathcal{L}f(x) = \int [f(x + F(x, z)) - f(x) - f'(x)F(x, z)]\lambda(dz), \tag{3.9}$$

we then see that

$$f(X_t) - f(X_0) - \int_0^t \mathcal{L}f(X_{s-})ds$$

is a martingale. Since $X$ has only countably many jumps, then the Lebesgue measure of the set of times where $X_s \ne X_{s-}$ is 0, and hence $f(X_t) - f(X_0) - \int_0^t \mathcal{L}f(X_s)ds = f(X_t) - f(X_0) - \int_0^t \mathcal{L}f(X_{s-})ds$. Therefore

$$f(X_t) - f(X_0) - \int_0^t \mathcal{L}f(X_s)ds$$

is a martingale.

These calculations (and analogous ones when $\sigma$ or $b$ in (3.7) is nonzero) are the motivation for what is known as the martingale problem. In this formulation we take $X_t$ to be the canonical coordinate process. One defines an integral operator by a formula such as (3.9) and then says that a probability measure $\mathbb{P}$ solves the martingale problem started at a point $x_0$ if $\mathbb{P}(X_0 = x_0) = 1$ and $f(X_t) - f(X_0) - \int_0^t \mathcal{L}f(X_s)ds$ is a martingale under $\mathbb{P}$ whenever $f \in C^2$. Note here that a solution is a probability, not a process.

Since the number of large jumps is finite in number over any finite time interval and the large jumps do not affect the existence or uniqueness of solutions, and one wants not to worry about the integrability of $\int |z|F(x, z)\lambda(dz)$, one often sees operators of the form

$$\mathcal{L}f(x) = \int [f(x + z) - f(x) - 1_{(|z|\le 1)}f'(x)z]n(x, dz). \tag{3.10}$$

If $n(x, dz)$ does not depend on $x$, the reader will note that in this case (3.10) is the infinitesimal generator of a Lévy process.

An example of an operator such as (3.10) is to let

$$n(x, dz) = \frac{a(x, z)}{|z|^{1+\alpha}} dz, \tag{3.11}$$

where $a$ is bounded above and below by positive constants, and $0 < \alpha < 2$. If $a$ is constant, we have the infinitesimal generator of a symmetric stable process of index $\alpha$.

The extension of the notion of martingale problem to higher dimensions causes no problems. One also often sees operators of the form $\mathcal{L} + \mathcal{M}$, where $\mathcal{L}$ is as in (3.10) and

$$\mathcal{M} f(x) = \tfrac{1}{2}\sigma(x)^2 f''(x) + b(x) f'(x).$$

This corresponds to adding a diffusion and drift term to the jump operator. Again, the higher dimensional analogues are what one expects.

*Pseudodifferential operators.* With $\mathcal{L}$ as in (3.10), let us see what we get when we apply $\mathcal{L}$ to the function $f(x) = e^{iux}$. We obtain

$$\mathcal{L} f(x) = e^{iux} \int [e^{iuz} - 1 - iuz 1_{(|z| \le 1)}] n(x, dz).$$

In the case where $n(x, dz)$ does not depend on $x$, i.e., $n(x, dz) = n(dz)$ for all $x$, this implies that the Fourier transform of $\mathcal{L}$ is $\psi(u)$, where

$$\psi(u) = \int [e^{iuz} - 1 - iuz 1_{(|z| \le 1)}] n(dz).$$

More generally, when $n(x, dz)$ does depend on $x$, we let

$$\psi(x, u) = \int [e^{iuz} - 1 - iuz 1_{(|z| \le 1)}] n(x, dz),$$

and we call $\psi(x, u)$ the symbol corresponding to the operator $\mathcal{L}$.

By the uniqueness of the Fourier transform, we can specify an operator by presenting $n(x, dz)$ or we can specify the operator by giving its symbol $\psi(x, u)$.

*Dirichlet forms.* Suppose one considers the operator $\mathcal{L} f(x) = (a(x) f'(x))'$. For this to make sense in terms of the usual notion of derivative, we need $a$ to be differentiable. If $g$ is $C^\infty$ with compact support, an integration by parts gives

$$- \int g(x) \mathcal{L} f(x)\, dx = \int g'(x) a(x) f'(x)\, dx.$$

The expression on the right makes sense for any measurable and bounded $a$, as long as $f$ and $g$ are differentiable with compact support, say, and is called a Dirichlet form, written $\mathcal{E}(f, g)$.

Under minimal hypotheses, a Dirichlet form will determine the process. See [FOT] for the basic theory of Dirichlet forms.

If we define

$$\mathcal{E}(f,f) = \int\int [f(y)-f(x)]^2 J(dx,dy), \tag{3.12}$$

where the measure $J$ is symmetric, the associated process will be of jump type. A special case is when $J(dx\,dy) = |x-y|^{-1-\alpha}dx\,dy$, and the associated process is a symmetric stable processes of index $\alpha$. We will briefly consider some properties of processes associated to the Dirichlet form given in (3.12), called symmetric jump processes, in the case where

$$J(dx,dy) = \frac{a(x,y)}{|x-y|^{1+\alpha}} dx\,dy \tag{3.13}$$

and also higher dimensional analogues, where $0 < \alpha < 2$, $a$ is bounded above and below by positive constants, and $a(x,y) = a(y,x)$ for all $x$ and $y$. Here the questions concern not uniqueness, since that is covered by the general theory of Dirichlet forms, but instead properties of the corresponding process. (But see also [AS].)

## 4. Pathwise uniqueness

When is the solution to (3.2) or (3.7) pathwise unique? When the coefficients are Lipschitz coefficients, then the standard Picard iteration procedure proves uniqueness. For example, for (3.7) we have the following, due to [Sk].

**Theorem 4.1** *If $\sigma$ and $b$ are bounded and Lipschitz, $\int \sup_x |F(x,z)|^2\lambda(dz) < \infty$ and*

$$\int |F(x,z)-F(y,z)|^2\lambda(dz) \le c_1|x-y|^2$$

*for all $x, y$, then there exists a solution to (3.7) and that solution is pathwise unique.*

Concerning (3.2), again, if $a$ is bounded and Lipschitz continuous, one would expect pathwise uniqueness, and indeed that is the case. In one dimension, however, in view of the result of Yamada-Watanabe [YW] for SDEs with respect to a one dimensional Brownian motion, one would hope that much weaker conditions would suffice for uniqueness. For example, the Yamada-Watanabe sufficient condition for pathwise uniqueness for diffusions is that $\int_0^\varepsilon \rho(x)^{-2}dx = \infty$ for all $\varepsilon > 0$, where $\rho$ is the modulus of continuity: $|\sigma(x)-\sigma(y)| \le \rho(|x-y|)$ for all $x$ and $y$.

For solutions to (3.2) we have the following; see [Km2] and [Ba4].

**Theorem 4.2** *Suppose $\alpha > 1$, $Z_t$ is a one dimensional symmetric stable process of index $\alpha$, and $|a(x)-a(y)| \le \rho(|x-y|)$ for all $x$ and $y$ and*

$$\int_0^\varepsilon \frac{1}{\rho(x)^\alpha} dx = \infty$$

*for every $\varepsilon > 0$. Then for every $x_0 \in \mathbb{R}$, the solution to*

$$dX_t = a(X_{t-})dZ_t, \qquad X_0 = x_0,$$

*is pathwise unique.*

This result is sharp. As a corollary, if $a$ is Hölder continuous of order $1/\alpha$, then we have pathwise uniqueness. Note that the smaller $\alpha$, the more continuity is required.

What happens when $\alpha \le 1$? When $\alpha < 1$, the paths of $Z_t$ are of bounded variation on finite time intervals, and I initially thought only continuity and appropriate boundedness of $a$ sufficed [Ba4]. However the solution constructed there, although measurable with respect to the $\sigma$-field generated by $Z$, is not adapted. In fact, when $\alpha \le 1$, the condition that $a$ be Lipschitz continuous turns out to be sharp; see [BBC]. Probably the Lipschitz continuity condition can be weakened by a logarithm term, but it does not suffice to let $a$ be Hölder continuous of order $\beta$ for any $\beta < 1$.

For diffusions Nakao [Na] showed that bounded variation of $\sigma$ suffices for (3.1) to have a unique solution, and LeGall [LG] improved this to $\sigma$ having finite quadratic variation. What is the appropriate analogue of this for the equation (3.2)?

A paper concerning pathwise existence along somewhat different lines is Williams [Wi]. Here the stochastic integral is not the Itô integral, but instead the rough integral of Lyons (cf. [Ly]), and the issue is existence.

Other papers related to pathwise properties of SDEs with jumps are [Ro] on backward SDEs, [JMW] on approximations of solutions, and [AT], [Fu], [FK], and [Ku] on flows of solutions.

## 5. Martingale problems

In this section let us discuss existence and solutions to martingale problems such as those mentioned in Section 3. When there exists a non-degenerate diffusion component, that term dominates the situation, and uniqueness holds. We have the following theorem of Stroock [St]

**Theorem 5.1** *Suppose*

$$\mathcal{L}f(x) = \sum_{i,j=1}^d a_{ij}(x)\frac{\partial^2 f}{\partial x_i \partial x_j}(x) + \sum_{i=1}^d b_i(x)\frac{\partial f}{\partial x_i}(x) \tag{5.1}$$

$$+ \int_{\mathbb{R}^d \setminus \{0\}} [f(x+z) - f(x) - 1_{(|z|\le 1)} \nabla f(x) \cdot z] n(x, dz) \tag{5.2}$$

*and*

*(a) $a$ is bounded and continuous and strictly elliptic;*

*(b) $b$ is bounded and measurable;*

*(c) $\int_A \frac{|z|^2}{1+|z|^2} n(x,dz)$ is bounded and continuous for each $A \subset \mathbb{R}^d \setminus \{0\}$.*

*Then there is a unique solution to the martingale problem for $\mathcal{L}$ started at any $x_0 \in \mathbb{R}^d$.*

This result should also be compared with that of Komatsu [Km1] and that of Lepeltier and Marchal [LM].

Suppose now that there is no diffusion component present. For early works see [Ts1], [Ts2], and [TTW]. A special case of a result of Komatsu [Km4] is the following.

**Theorem 5.2** *Suppose*

$$\mathcal{L}f(x) = \int_{\mathbb{R}^d \setminus \{0\}} [f(x+z) - f(x) - 1_{(|z|\le 1)} \nabla f(x) \cdot z] \Big[ \frac{1}{|z|^{d+\alpha}} dz + n_0(x, dz) \Big]$$

*for some* $\alpha \in (0,2)$*,* $n^*(dz) = \sup_x n(x, dz)$*, and*

$$\int (1 \wedge |z|^\alpha) n^*(dz) < \infty.$$

*Then there is a unique solution to the martingale problem for* $\mathcal{L}$ *started at any point* $x_0$*.*

See also [NT].

Martingale problems are closely related to the notion of weak uniqueness for SDEs. Concerning weak existence of solutions for (3.2) see [PZ], [Za1], and [Za2].

Another approach that has been explored has been the use of pseudodifferential operators. See works of Komatsu, Jacob, Hoh, and Schilling: [Km3], [HJ1], [HJ2], [Ho1], [Ho2], [Ho3], [Ho4], [JL], [Ja1], [Ja2], [Ja3], [JS], [Ja4], and [Ja5].

A recent result along these lines is that of Kolokoltsov [Kl2]. If $\psi_i(u)$, $i = 1, \dots, n$, are the symbols associated to Lévy processes satisfying certain conditions and $a_i(x)$ are suitable functions of $x$, then there is a unique solution to the martingale problem associated to the symbol

$$\psi(x, u) = \sum_{i=1}^n a_i(x) \psi_i(u).$$

Can one obtain a similar result when we have an infinite sum instead of a finite one?

In most of the above results, the operators are the sum of a finite number of integral terms plus a lower order term. There are fewer papers that handle variable order terms without assuming a considerable amount of continuity in the $x$ variable. One such paper is [Ba1]. This paper deals with one-dimensional processes, and the exact conditions are rather complicated. One can give an example, though, that illustrates the theorem. Suppose $\mathcal{L}$ is given by (3.10), where

$$n(x, dz) = \frac{1}{|z|^{1+\alpha(x)}} dz. \tag{5.3}$$

Qualitatively, this says that at a point $x$, the process behaves like a symmetric stable process of order $\alpha(x)$. Note that the order depends on $x$. For this example, the conditions of [Ba1] reduce to **Theorem 5.3** *Suppose $\alpha(x)$ is Dini continuous, bounded above by a constant less than 2 and bounded below by a constant greater than 0. If $\mathcal{L}$ is given by (3.10) with $n(x, dz)$ defined by (5.3), then there exists a unique solution to the martingale problem for $\mathcal{L}$ started at any point $x_0 \in \mathbb{R}$.*

Dini continuity of $\alpha$ means that there exists $\rho$ such that $|\alpha(x) - \alpha(y)| \leq \rho(|x - y|)$ for all $x$ and $y$ and

$$\int_0^\varepsilon \frac{\rho(x)}{x} dx < \infty$$

for all $\varepsilon > 0$.

Other papers that deal with properties of processes associated to operators of variable order without assuming a great deal of smoothness include [Ba2], [Ne], [Ts3], [Km5], [Kl1], [Ue].

It would be greatly desirable to have a uniqueness theorem for variable order operators where the hypotheses are relatively simple.

## 6. Regularity of harmonic functions

In this section we give some regularity results for harmonic functions related to the integral operators $\mathcal{L}$ such as those similar to (3.11). Before doing this, we mention that there has been a great deal of recent work on the potential theory of symmetric stable processes in domains, such as Fatou theorems, boundary Harnack principle, intrinsic ultracontractivity, Green function estimates, etc. See [Ch] for a survey of some of the recent work. Little of this has been extended to more general jump processes, but it would be worthwhile doing so.

First we consider the Harnack inequality. Suppose $\mathcal{L}h(x) = 0$ in a domain $D \subset \mathbb{R}^d$. Here the operators $\mathcal{L}$ that we are considering are of the form

$$\mathcal{L}f(x) = \int_{\mathbb{R}^d \setminus \{0\}} [f(x+z) - f(x) - 1_{(|z| \leq 1)} \nabla f(x) \cdot z] n(x, z) dz, \qquad (6.1)$$

where

$$n(x, z) = \frac{a(x, z)}{|z|^{d+\alpha}}, \qquad (6.2)$$

$\alpha \in (0, 2)$, and $a$ is bounded above and below by positive constants. It may so happen that what one would like to consider as a harmonic function is not regular enough to be in the domain of $\mathcal{L}$. So a more general definition of harmonic function is to say that $h$ is harmonic if $h(X_{t \wedge \tau_D})$ is a martingale. Here $\tau_D = \inf\{t : X_t \notin D\}$, the first exit time of the domain $D$, and $X_t$ is the Markov process associated to the operator $\mathcal{L}$. We have the following Harnack inequality [BL1].

**Theorem 6.1** *Suppose $\mathcal{L}$ is given by (6.1), $n$ by (6.2), $n(x, z) = n(x, -z)$ for all $x$ and $z$, $x_0 \in \mathbb{R}^d$, $R > 0$, $h$ is nonnegative and harmonic in a ball $B(x_0, 2R)$ of radius $2R$ centered at $x_0$, and $h$ is bounded in $\mathbb{R}^d$. Then there exists a constant $c_1$ not depending on $h$, $x_0$, or $R$ such that*

$$h(x) \leq c_1 h(y), \qquad x, y \in B(x_0, R).$$

This theorem has been extended by [SV] to more general $n$, but still ones that are essentially stable of some order $\alpha$.

One of the reasons one is interested in Harnack inequalities is that they often imply that harmonic functions must be continuous. In [BL1] it is also proved that

**Theorem 6.2** *Under the same conditions as Theorem 6.1 there exist constants $c_2$ and $\gamma \in (0,1)$ such that*

$$|h(x) - h(y)| \le c_2 \Big( \frac{|x-y|}{R} \Big)^\gamma \|h\|_\infty.$$

In other words, harmonic functions are Hölder continuous.

What if the function $n$ is not of the form (6.2)? Suppose there exist constants $c_3, c_4$ and $0 < \alpha < \beta < 2$ such that

$$\frac{c_3}{|x-y|^{d+\alpha}} \le n(x,y) \le \frac{c_4}{|x-y|^{d+\beta}}, \qquad |x-y| \le 1, \tag{6.3}$$

and with some appropriate conditions on $n(x,y)$ when $|x-y| \ge 1$? In [BK] it was proved that the conclusion of Theorem 6.1 holds provided that $\beta - \alpha < 1$ and one allows $c_1$ to depend on $R$. It was shown that the dependence on $R$ cannot be dispensed with.

It is still unknown whether a Harnack inequality must hold in the case $\beta - \alpha \ge 1$, or whether a counterexample exists.

If one turns to symmetric jump processes, that is, ones determined by the Dirichlet form (3.12), that the analogue of Theorem 6.1 holds is implicitly proved in [BL2]; see also [CK]. In the paper [BL2] upper and lower bounds on the transition probability densities were found, extending and improving results of Komatsu [Km5].

I believe there is also an analogue of [BK], but that in this case the full range $0 < \alpha < \beta < 2$ is allowed. Whether the analogue of Theorem 6.2 holds in this full range is at present unknown.

If one assumes more regularity in the operator $\mathcal{L}$, that is, stronger conditions on $n$, then more can be said about regularity. See [Kl2].

## Acknowledgment

I would like to thank Jim Pitman for his encouragement to write this paper, Moritz Kassmann for his very helpful suggestions, and David Applebaum for graciously providing me with material from his book.